\documentclass[a4paper,10pt]{article}

\usepackage[T1]{fontenc}
\usepackage[latin1]{inputenc}
\usepackage[english]{babel}
\usepackage{multicol}
\usepackage{fancyhdr}
\usepackage{amsmath,amssymb,amsthm}
\usepackage[all]{xy}
\usepackage{stackrel}
\usepackage{geometry}

\renewcommand\c[1]{\mathcal{#1}}
\newcommand\db{\c D^b}
\def\t#1#2{^{#1}\widetilde{T}_{#2}}

\title{Simple connectedness of quasitilted algebras}

\author{Patrick Le Meur
\footnote{\textit{adress:} Département de Mathématiques, \'Ecole normale supérieure de
  Cachan, 61 avenue du Président Wilson, 94235 CACHAN cedex, FRANCE}
\footnote{\textit{e-mail:} plemeur@dptmaths.ens-cachan.fr}
}

\date{\today}

\newtheorem{prop}{Proposition}[section]
\newtheorem{Thm}[]{Theorem}
\newtheorem{lem}[prop]{Lemma}

\newtheorem{rem}[prop]{Remark}

\begin{document}
\maketitle
\abstract{
Let $A$ be a basic connected finite dimensional algebra over an algebraically closed
field. Assuming that $A$ is quasitilted, we prove that $A$ is simply connected if and only if
$HH^1(A)=0$. This generalises a result of I.~Assem, F.~U.~Coelho and
S.~Trepode and which proves the same equivalence for tame quasitilted
algebras.
}

\section*{Introduction}
Let $A$ be a finite dimensional algebra over an algebraically closed
field $k$. In order to study the category $mod(A)$ of finite
dimensional (right)
$A$-modules we may assume that $A$ is basic and connected. In this
study, covering techniques introduced in \cite{bongartz_gabriel} and
\cite{riedtmann} have proved to be a very powerful tool. Indeed, a
Galois covering $\c C\to A$ (with $\c C$ a locally bounded
$k$-category) reduces the study of part of
$mod(A)$ to the one of $mod(\c C)$
which is easier to handle (see for example \cite{gabriel}). From this
point of view, simply connected algebras are of particular
interest. Recall that if $Q$ is the ordinary quiver of $A$ and if $kQ$
is the path algebra of $Q$, then there exists a (non necessarily
unique) surjective algebra
morphism (or presentation) $\nu\colon kQ\twoheadrightarrow A$ (see \cite{ars} for
example). Moreover, given such a presentation, one can define the
fundamental group $\pi_1(Q,Ker(\nu))$ of $\nu$ (see
\cite{martinezvilla_delapena}). With this setting, the algebra $A$ is
called simply connected if and only if $Q$ has not oriented cycle
(\textit{i.e.} no non trivial oriented path whose source equals its
target, the algebra is then called triangular) and $\pi_1(Q,Ker(\nu))=1$
for any $\nu\colon kQ\twoheadrightarrow 
A$ (see \cite{assem_skowronski}). Equivalently (\cite{skowronski} and
\cite{lemeur2}) $A$ is simply connected if and only if there 
exists no Galois covering $\c C\to A$ with non trivial group and with
$\c C$ a connected locally bounded $k$-category.

To prove that $A$ is simply connected seems to be a difficult problem,
\textit{a priori}, since one has to check that various groups are
trivial. Hence, it is worth looking for a simpler characterisation of
 simple connectedness. It was asked by A.~Skowro\'nski
(\cite{skowronski2}) whether the equivalence ``$A$ is simply
connected if and only if  $HH^1(A)=0$'' is satisfied for $A$ a tame
triangular algebra. This equivalence is true for
tilted algebras (see \cite{assem_marcos_delapena} for the tame case
and \cite{lemeur6} for the general case), for piecewise hereditary
algebras of type any quiver (see \cite{lemeur6}), for tame
quasititled algebras (see \cite{assem_coelho_trepode}) and it is
conjectured (\textit{loc.cit.}) that this equivalence is
true for any quasitilted algebra. 

Recall (\cite{happel_reiten_smalo}) that a quasitilted algebra is an
algebra isomorphic to $End_{\c H}(T)^{op}$ where $\c H$ is a
hereditary abelian $k$-linear category (with finite dimensional Hom
and Ext spaces) and where $T\in\c H$ is a (basic) tilting
object. In particular, a quasitilted algebra has global dimension at
most $2$ (see \textit{loc.cit.}). Quasitilted algebras were introduced in order to give a common
framework to the class of tilted algebras (introduced in
\cite{happel_ringel}) and to the class of canonical algebras
(introduced in \cite{ringel2}).
In this text, we prove the following result:
\begin{Thm}
\label{thm1}
Let $A$ be a basic connected finite dimensional $k$-algebra. If $A$ is
quasitilted, then:
\begin{equation}
  \text{$A$ is simply connected}\ \Leftrightarrow\ HH^1(A)=0\notag
\end{equation}
Moreover, if $A$ is tilted of type $Q$, then $A$ is simply connected
if and only if $Q$ is a tree.
\end{Thm}
Hence, the above theorem solves the above conjecture of
\cite{assem_coelho_trepode} and it also answers positively the
above question of A.~Skowro\'nski (\cite{skowronski2}) for quasitilted
algebras (of finite, tame or wild type). Recall that in Theorem~\ref{thm1}, the case of tilted
algebras and the one of quasitilted algebras which are derived equivalent to a
hereditary algebra have been successfully treated in \cite{assem_marcos_delapena} and
\cite{lemeur6}. Here, we say hat two algebras are derived equivalent
if and only if their derived categories of bounded complexes of finite
dimensional modules are triangle equivalent.

In order to prove Theorem~\ref{thm1}, we use ideas from
\cite{lemeur6}. More precisely, given a quasitilted algebra $A$ which
is not derived equivalent to a  hereditary algebra, we find a suitable
algebra $B$ which is derived equivalent to $A$ and for which the equivalence of
Theorem~\ref{thm1} may be proved easily. Then, we prove that $A$
is simply connected if and only if $B$ is simply connected by
establishing a correspondence between the Galois coverings of $A$ and
those of $B$. This correspondence is very
similar to the one of \cite{lemeur6} (and of \cite{lemeur5}) since we
shall compare the Galois coverings of $A$ and those of
$End_{\db(A)}(T)$ for some $T\in\db(A)$ (where $\db(A)$ denotes the
derived category of bounded complexes of $A$-modules). Recall that, in
\cite{lemeur6}, the suitable algebra $B$ associated a
piecewise hereditary algebra $A$ of type $Q$ was chosen to be the path
algebra $kQ$. Here, we shall take for $B$ a squid algebra (see
\cite{ringel}). Indeed, it
was proved in \cite{happel_reiten} that a hereditary abelian
$k$-linear category with tilting object  and which is not derived equivalent to the module
category of a hereditary algebra is derived equivalent to a squid algebra.

The text is organised as follows. In Section~$1$ we fix some
notations. In Section~$2$, we construct the
above correspondence. In Section~$3$, we prove the
Theorem~\ref{thm1} for squid algebras. Finally, Section~$4$ is devoted
to the proof of this theorem.

\section{Notations}

A $k$-category $\c C$ is a category whose collection $ob(\c C)$ of
objects is a set, whose space of morphisms $_y\c C_x$ (or $\c C(x,y)$)
from $x$ to $y$ is a $k$-vector space for any $x,y\in ob(\c C)$, and whose composition of
morphisms is $k$-bilinear. All functors between $k$-categories will be
assumed to be $k$-linear.
A basic connected finite dimensional $k$-algebra $A$ will always be
considered as a locally bounded $k$-category (see
\cite{bongartz_gabriel}) with set of objects a complete set
$\{e_1,\ldots,e_n\}$ of pairwise orthogonal primitive idempotents,
with space of morphisms from $e_i$ to $e_j$ equal to $e_jAe_i$ and with
composition of morphisms induced by the product of $A$.\\

Following \cite{bongartz_gabriel}, a (right) module over a locally
bounded $k$-category $\c C$ is a $k$-linear covariant functor from $\c C^{op}$ to
the category of $k$-vector spaces. Such a module $M$ is called finite
dimensional if $\sum_{x\in ob(\c C)}dim_k\ M(x)<\infty$. In
particular, for any $x\in ob(\c C)$, the indecomposable projective
module $y\mapsto\ _x\c C_y$ will be denoted by $_x\c C_?$. The category
of finite dimensional $\c C$-modules will be denoted by $mod(\c
C)$. The derived category of bounded
complexes of $\c C$-modules will be denoted by $\db(\c
C)$ and $\Sigma$ will denote the shift functor. Recall that if $\c C$ has finite
global dimension, then $\db(\c C)$ is
equivalent to the homotopy category of bounded complexes of finite
dimensional projective
$\c C$-modules. The Auslander-Reiten translation (see
\cite{happel_book}) on $\db(\c C)$ will be denoted by $\tau_{\c
  C}$. Also, if $\c H$ is a hereditary abelian category with tilting
objects, then  we shall
write $\tau_{\c H}$ for the Auslander-Reiten translation on $\c H$ and
on $\db(\c H)$.\\

For a reminder on Galois coverings, we refer the reader to
\cite{bongartz_gabriel} or \cite{cibils_marcos}. A Galois covering
$F\colon \c C\to A$ will be called connected if and only if $\c C$
(and therefore $A$) is connected and locally bounded. Recall that if
$F\colon\c C\to A$ is a Galois covering with group $G$ and with $\c C$
and $A$ locally bounded, then $F$ defines a triangle functor
$F_{\lambda}\colon\db(\c C)\to\db(A)$ (see for example
\cite[Lem. 2.1]{lemeur6}). Moreover, the group $G$ acts on
$\db(\c C)$ by triangle isomorphisms $(g,X)\in G\times \db(\c
C)\mapsto\ ^gX$. For this action, $F_{\lambda}$ is $G$-invariant and
for any $X,Y\in \db(\c C)$, the following maps induced by
$F_{\lambda}$ are linear isomorphisms:
\begin{equation}
  \bigoplus\limits_{g\in G}Hom_{\db(\c C)}(\,^gX,Y)\xrightarrow{\sim}
  Hom_{\db(A)}(F_{\lambda}X,F_{\lambda}Y),\ \ \bigoplus\limits_{g\in G}Hom_{\db(\c C)}(X,\,^gY)\xrightarrow{\sim}
  Hom_{\db(A)}(F_{\lambda}X,F_{\lambda}Y)\notag
\end{equation}
For short, these properties on $F_{\lambda}$ will be called the
\textit{covering properties of $F$}. Recall (\cite[Lem. 4.1]{lemeur6}) that
$F_{\lambda}$ verifies $\tau_A\circ F_{\lambda}\simeq F_{\lambda}\circ
\tau_{\c C}$. An indecomposable object $X\in\db(A)$ is called of the
first kind w.r.t. $F$ if and only if $X\simeq
F_{\lambda}\widetilde{X}$ for some $\widetilde{X}\in\db(\c C)$ (which
is necessarily indecomposable). More
generally, an object $X\in\db(A)$ is called of the first kind w.r.t. $F$ if and
only if $X$ is the direct sum of indecomposable objects of the first
kind w.r.t. $F$. Finally, for $T\in\db(A)$, we introduce two
assertions depending on $T$ and $F$ and which will be used in this
text:

\begin{enumerate}
\item[($H_1$)] $T$ is of the first kind w.r.t. $F$.
\item[($H_2$)] for every indecomposable direct summand $X\in\db(A)$ of
  $T$,  for any $\widetilde{X}\in\db(\c C)$ such that
  $F_{\lambda}\widetilde{X}\simeq X$ in $\db(A)$, and for any  
    $g\in G\backslash\{1\}$, we have
    $^g\widetilde{X}\not\simeq\widetilde{X}$ in $\db(\c C)$.
\end{enumerate}

For a reminder on tilting objects in hereditary abelian categories, we
refer the reader to \cite{happel_reiten_smalo}, for a reminder on
cluster categories and on cluster tilting objects we refer the reader
to \cite{bmrrt}. The cluster category of a finite dimensional algebra
$A$ (resp. of a hereditary abelian category $\c H$) will be denoted by
$C_A$ (resp. $C_{\c H}$).\\

If $\c T$ is a triangulated category with shift functor
$\Sigma$, we set $Ext^i_{\c T}(?,!):=Hom_{\c T}(?,\Sigma^i!)$. Also,
if $A$ is a finite dimensional algebra, we shall write $Ext^i_A$
instead of $Ext^i_{\db(A)}$, for simplicity.\\

Finally, if $\c A$ is a dg category, $Dif\ \c A$ will denote the dg
category of dg $\c A$-modules, $\c D(\c A)$ will denote the associated
derived category. Recall that if $\c A$ is a $k$-category considered
as a dg category concentrated in degree $0$, then $\c D(\c A)$ is the
usual derived category of (unbounded) complexes of $\c A$-modules. If
$\c B$ is another dg category and if $X$ is a dg
$\c B-\c A$-bimodule, $?\stackrel[\c B]{}{\otimes}X\colon Dif\ \c B\to
Dif\ \c A$ will denote the tensor product dg functor and
$?\stackrel[\c B]{\mathbb{L}}{\otimes}X\colon \c D(\c B)\to
\c D(\c A)$ will denote its left derived functor. For a reminder on dg
categories, we refer the reader to \cite{keller_dg}.

\section{Invariance of simple connectedness under tilting}
Let $A$ be a basic connected finite dimensional $k$-algebra. Let
$\db(\c H)\xrightarrow{\sim}\db(A)$ be a triangle equivalence
where $\c H$ is a hereditary abelian category with tilting objects.  Finally, let $T\in\db(A)$ be a basic object
 such that:
 \begin{enumerate}
 \item $T$ is a cluster tilting object of
$C_A$
\item $Ext_A^i(T,T)=0$ for any $i\neq 0$.
 \end{enumerate}
  In this section, we shall compare the Galois coverings of $A$ and
  those of $A':=End_{\db(A)}(T)$ in order to prove the following implication:
\begin{equation}
  \text{$A'$ is simply connected}\ \Rightarrow\
  \text{$A$ is simply connected}\tag{$\star\star$}
\end{equation}
Set $T=T_1\bigoplus\ldots\bigoplus T_n\in \db(A)$ with
$T_1,\ldots,T_n\in \db(A)$ indecomposables (where $n=rk(K_0(A)$). Recall that $A'$ is a
locally bounded $k$-category with set of objects $\{T_1,\ldots,T_n\}$,
with space of morphisms from $T_i$ to $T_j$ equal to $Hom_{\db(A)}(T_i,T_j)$
and with composition of morphisms induced by the composition in
$\db(A)$.

\subsection{(Cluster) tilting objects of the first kind w.r.t. Galois coverings}

In order to prove ($\star\star$), we shall
associate Galois coverings of $A'$ to Galois coverings of $A$ using a
construction of \cite[Sect. 2]{lemeur5} and then use the characterisation
\cite[Cor. 4.5]{lemeur2} of simple connectedness in terms of Galois coverings. In this purpose, the
following lemma will be useful. Its proof is based on the work made in \cite{lemeur6}.
\begin{lem}
\label{lem1}
 (see \cite[Prop. 6.5, Prop. 6.8]{lemeur6})
  Let $F\colon \c C\to A$ be a Galois covering with group $G$ and with
  $\c C$ locally bounded. Then, ($H_1$) and ($H_2$) are satisfied for
  $F$ and for any object of $\db(A)$ which is a cluster tilting object
  of $C_A$.
\end{lem}
\noindent{\textbf{Proof:}} For simplicity, we shall make no
distinction between an object and its isomorphism class. Let $\c
S\subseteq \db(A)$ be the class of objects $R\in\db(A)$ which are
isomorphic to the image of a
cluster tilting object of $C_{\c H}$ under the equivalence $\db(\c
H)\to \db(A)$. Hence, $\c S$ is the
  class of cluster tilting objects of $C_A$. In particular, it contains
  $A$.  Let $\sim$ be the equivalence relation on $\c S$
  generated by the following property: ``\textit{if $R,R'\in\c S$ are
    such that $R=X\bigoplus\overline{R}$, $R'=Y\bigoplus\overline{R}$
    with $X,Y\in\db(A)$ indecomposables and verifying at least one of
    the following properties:
    \begin{enumerate}
    \item $X\simeq (\tau_A\Sigma^{-1})^mY$ for some $m\in\mathbb{Z}$,
    \item there exists a triangle $X\to M\to Y\to\Sigma X$ of $\db(A)$
      with $M\in
      add(\overline{R})$,
    \item there exists a triangle $Y\to M\to X\to\Sigma Y$ of $\db(A)$
      with $M\in
      add(\overline{R})$.
    \end{enumerate}
then $R\sim R'$''
}. Here, $add(R)$ denotes the full additive subcategory of $\db(A)$
closed under isomorphisms and generated by the indecomposable direct
summands of $R$. Since any cluster tilting object of $C_{\c H}$ is isomorphic (in
$C_{\c H}$) to a tilting object of $\c H$ (see
\cite[Sect. 3]{bmrrt}), since the Hasse diagram of tilting objects
of $\c H$ is connected (see \cite{happel_unger3}), and since
$\db(\c H)\xrightarrow{\sim}\db(A)$ commutes with $\Sigma$
and preserves Auslander-Reiten triangles, we infer that $\c
S$ is an equivalence class for $\sim$. On the other hand,
$F_{\lambda}\colon\db(\c C)\to\db(A)$ commutes with $\Sigma$ and is
compatible with $\tau_{\c C}$ and $\tau_A$ (\textit{i.e.} $F_{\lambda}\circ \tau_{\c
  C}\simeq \tau_A\circ F_{\lambda}$, see
\cite[Lem. 4.1]{lemeur6}). Therefore, using \cite[Prop 6.5]{lemeur6}
and \cite[Prop 6.8]{lemeur6}, we deduce that if $R\sim R'$, then the
conclusion of the lemma holds for $R$ and $F$ if and only if it holds for
$R'$ and $F$. Since ($H_1$) and ($H_2$) are clearly satisfied for $T=A$, the lemma is proved.
\hfill$\square$\\

 \subsection{Galois coverings of $A'$ associated to Galois coverings
   of $A$}

Now we can recall the construction of \cite[Sect. 2]{lemeur5} which
associates Galois coverings of $A'$ to Galois coverings of $A$. Fix
$F\colon \c C\to A$ a Galois covering with group $G$ and with $\c C$
locally bounded. Assume that there exist
$\widetilde{T}_1,\ldots,\widetilde{T}_n\in\db(\c C)$ together with
isomorphisms $\lambda_i\colon
F_{\lambda}\widetilde{T}_i\xrightarrow{\sim}T_i$ in $\db(A)$, for
every $i$ (see Lemma~\ref{lem1}). Then, we define $\c C'$ to be the following
  $k$-category:
  \begin{enumerate}
  \item the set of objects of $\c C'$ is $\{\ ^g\widetilde{T}_i\ |\
    g\in G,\ i\in\{1,\ldots,n\}\}$.
\item $_{^h\widetilde{T}_j}\c C'_{^g\widetilde{T}_i}:= Hom_{\db(\c
    C)}(\,^g\widetilde{T}_i,\,^h\widetilde{T}_j)$ for any $g,h\in G$
  and $i,j\in\{1,\ldots,n\}$.
\item the composition in $\c C'$ is induced by the composition in
  $\db(\c C)$.
  \end{enumerate}
Hence, $\c C'$ is the full subcategory of $\db(\c C)$ whose objects
are the complexes $\t g i$. Moreover, we define a $k$-linear functor
$F'\colon \c C'\to A'$ as
follows:
\begin{equation}
  \begin{array}{crcl}
    F'\colon & \c C'&\to & A'\\
    & ^g\widetilde{T}_i\in ob(\c C') & \mapsto & T_i\in ob(A')\\
    & u\in\ _{^h\widetilde{T}_j}\c C'_{^g\widetilde{T}_i} & \mapsto &
    T_i\xrightarrow{\lambda_j\circ F_{\lambda}u\circ\lambda_i^{-1}} T_j
  \end{array}\notag
\end{equation}
The following lemma was proved in \cite{lemeur5} in the case $T\in
mod(A)$. However, the reader may easily check that the proof still
works in our situation ($T\in\db(A)$):
\begin{lem}
\label{lem2}
  (see \cite[Rem. 2.1, Lem. 2.2]{lemeur5})
The $G$-action on $\db(\c C)$ naturally defines a $G$-action on $\c
C'$. For this action, $F'\colon \c C'\to A'$ is a Galois covering with group $G$ and $\c C'$
is a locally bounded $k$-category.
\end{lem}

\begin{rem}
  \label{rem1}
Since $Ext_A^m(T,T)=0$ for any $m\neq 0$ and since $F_{\lambda}$ has
the covering property, we infer that $Ext_{\c C}^m(\,\t g i,\,\t h
j)=0$ for any $g,h\in G$, $i,j\in\{1,\ldots,n\}$ and  $m\in\mathbb{Z}\backslash\{0\}$.
\end{rem}

\subsection{Connectedness of Galois coverings}

Let us keep the notations of the preceding subsection. Since we are
interested in \textit{connected} Galois coverings, we
need to check when $\c C'$ is connected. In this purpose, we shall
prove the following proposition.
\begin{prop}
\label{prop1}
  $\c C$ and $\c C'$ are derived equivalent. In particular, $\c C'$ is connected if and only if $\c C$ is connected.
\end{prop}

We shall prove Proposition~\ref{prop1} in two steps: first 
we construct a fully faithful triangle functor $\Psi\colon\db(\c
C')\to\db(\c C)$ which maps the indecomposable projective $\c
C'$-module $_{\t g i}\c C'_?$ to an object of $\db(\c
C)$ isomorphic to $\t g i\in\db(\c C)$. Then, we prove that this
functor is dense. 
\begin{lem}
\label{lem3}
  There exists $\Psi\colon\db(\c C')\to\db(\c C)$ a fully faithful
  triangle functor such that $\Psi(\,_{\t g i}\c C'_?)\simeq\ \t g i$
  in $\db(\c C)$, for any $g\in G$, $i\in\{1,\ldots,n\}$. Moreover,
  $\Psi$ has a right adjoint triangle functor $\db(\c C)\to\db(\c C')$.
\end{lem}
\noindent{\textbf{Proof:}} 
We may assume that $\t g i=\t g i ^{\bullet}$ is a bounded complex of projective $\c
C$-modules, for any $g,i$.

$\bullet$ \textbf{A dg category $\c B$ derived equivalent to $\c C'$}.
Denote by $\c B$ the following dg category:
\begin{enumerate}
\item the set of objects is $\{\
\t g i \ |\ g\in G,\ i\in\{1,\ldots,n\}\}$,
\item $\c B^d(\,\t g i,\,\t h j):=\left\{(f_m\colon\ \t g i ^m\to \ \t h
    j ^{m+d})_{m\in\mathbb{Z}}\ |\ f_m\ \text{is a morphism of $\c
      C$-modules}\right\}$,
\item the differential $df$ of $f=(f_m)_{m\in\mathbb{Z}}\in\c B^d(\,\t g i,\,\t h j)$
  is given by:
  \begin{equation}
    (df)_m=d_{\t h j}^{m+d}\circ f_m-(-1)^d f_{m+1}\circ d_{\t g i}^m \notag
  \end{equation}
\end{enumerate}
 
Since $\t g i$ is a bounded complex of projective $\c C$-modules and
thanks to Remark~\ref{rem1},  there is an isomorphism of $k$-categories $H^0\c
B\xrightarrow{\sim}\c C'$ extending the identity map on
objects. To $\c B$ is associated the sub dg
category $\tau_{\leqslant 0}\c B$ with the same objects as $\c B$ and
such that $(\tau_{\leqslant 0}\c B)(X,Y)$ is truncated complex
$\tau_{\leqslant 0}(\c B(X,Y))$ for any $X,Y$. Thus, we have natural
dg functors:
\begin{equation}
  \c B\leftarrow \tau_{\leqslant 0}\c B\rightarrow H^0\c B\notag
\end{equation}
Once again, by assumption on $\t g i$ and thanks to Remark~\ref{rem1}, these functors induce isomorphisms of
graded categories:
\begin{equation}
  H^{\bullet}\c B\xleftarrow{\sim}H^{\bullet}\tau_{\leqslant 0}\c
  B\xrightarrow{\sim} H^0\c B\tag{$i$}
\end{equation}

On the other hand, $\tau_{\leqslant 0}\c B\to \c B$
(resp. $\tau_{\leqslant 0}\c B\to H^0\c B$) defines a dg $\tau_{\leqslant 0}\c
B-\c B$-bimodule $M$ (resp. a dg $
\tau_{\leqslant 0}\c B-H^0\c B$-bimodule $N$) such that $M(X,Y)=\c
B(X,Y)$ for any $X\in ob(\c B)$ and $Y\in  ob(\tau_{\leqslant 0}\c B)$
(resp. such that $N(X,Y)=H^0\c
B(X,Y)$ for any $X\in ob(H^0\c B)$ and $Y\in  ob(\tau_{\leqslant 0}\c
B)$). The bimodules $M$ and $N$ define triangle functors:
\begin{equation}
  \c D(\c B)\xleftarrow{?\stackrel[\tau_{\leqslant 0}\c B]{\mathbb{L}}{\otimes} M}\c
  D(\tau_{\leqslant 0}\c B)\xrightarrow{?\stackrel[\tau_{\leqslant 0}\c
    B]{\mathbb{L}}{\otimes} N}\c D(H^0\c B)\tag{$ii$}
\end{equation}

Using \cite[6.1]{keller_dg} and the isomorphisms of ($i$), we infer that
the above functors ($ii$) are triangle equivalences. Remark that since
$H^0\c B$ is concentrated in degree $0$, the derived category $\c
D(H^0\c B)$ is exactly the derived category of complexes of $H^0\c B$-modules.
Finally, for any $X\in ob(\c
B)=ob(\tau_{\leqslant 0}\c B)=ob(H^0(\c B))$, we have (\cite[6.1]{keller}):
\begin{enumerate}
\item $X^{\wedge}\stackrel[\tau_{\leqslant 0}\c B]{}{\otimes}M\simeq
  M(?,X)=X^{\wedge}$ in $Dif\ \c B$,
\item $X^{\wedge}\stackrel[\tau_{\leqslant 0}\c B]{}{\otimes}N\simeq
  N(?,X)=X^{\wedge}$ in $Dif\ H^0\c B$,
\end{enumerate}
so that:
\begin{enumerate}
\item $X^{\wedge}\stackrel[\c B]{\mathbb{L}}{\otimes} M\simeq
  X^{\wedge}$ in $\c D(\tau_{\leqslant 0}\c B)$,
\item $X^{\wedge}\stackrel[H^0\c
    B]{\mathbb{L}}{\otimes} N\simeq
  X^{\wedge}$ in $\c D(H^0\c B)$.
\end{enumerate}
These isomorphisms together with the equivalences ($ii$) and the isomorphism $H^0(\c
B)\simeq \c C'$ prove that there exists a triangle equivalence
$\Phi\colon \c
D(\c C')\xrightarrow{\sim} \c D(\c B)$ which maps $_{\t g i}\c C'_?$
to an object of $\c D(\c B)$ isomorphic to $\t g i
^\wedge$, for any $g,i$.

$\bullet$ \textbf{The triangle functor $?\stackrel[\c
B]{\mathbb{L}}{\otimes}\widetilde{T}\colon \c D(\c B)\to \c D(\c
C)$}. The complexes of $\c C$-modules $\t g i$ naturally define a dg $\c B-\c
C$-bimodule $\widetilde{T}$ such that $\widetilde{T}(x,\,\t g i)=\ \t g i
(x)$ for any $\t g i\in ob(\c B)$ and any $x\in ob(\c C)$. This
bimodule defines a triangle functor:
\begin{equation}
  ?\stackrel[\c B]{\mathbb{L}}{\otimes}\widetilde{T}\colon \c D(\c
  B)\to \c D(\c C)\notag
\end{equation}
Notice that \cite[6.1]{keller_dg} implies that:
\begin{equation}
  (\forall g,i)\ \ \t g i ^{\wedge}\stackrel[\c
B]{\mathbb{L}}{\otimes}\widetilde{T}\simeq \widetilde{T}(?,\,\t g i)=\ \t g i\tag{$iii$}
\end{equation}
Since $\t g i$ is a bounded complex of projective $\c C$-modules, we
infer (using \cite[6.2]{keller_dg}), that $?\stackrel[\c
B]{\mathbb{L}}{\otimes}\widetilde{T}$ admits a right adjoint triangle
functor $\c D(\c C)\to \c D(\c B)$. 

$\bullet$ \textbf{The triangle functor $\Psi\colon \c D^b(\c C')\to \c D^b(\c
C)$ and its right adjoint $\Theta\colon \c D^b(\c C)\to \c D^b(\c
C')$}.
Let us set $\Psi:=  ?\stackrel[\c
B]{\mathbb{L}}{\otimes}\widetilde{T}\circ \Phi\colon\c D(\c C')\to \c
D(\c C)$ and let us denote by $\Theta\colon \c D(\c C)\to \c D(\c C')$
the composition of a quasi inverse of $\Phi\colon\c D(\c
C')\xrightarrow{\sim}\c D(\c B)$ with the right adjoint $\c D(\c C)\to
\c D(\c B)$ of $?\stackrel[\c
B]{\mathbb{L}}{\otimes}\widetilde{T}$. Thus, the pair $(\Psi,\Theta)$ is
adjoint. Moreover, the construction of $\Phi$ and ($iii$) prove that:
\begin{equation}
  (\forall g,i)\ \ \Psi(\ _{\t g i}\c C'_?)\simeq\ \t g i\tag{$iv$}
\end{equation}
This proves that $\Psi$ maps $\db(\c C')$ into $\db(\c C)$ and that it
induces
a triangle functor $\Psi\colon\db(\c C')\to\db(\c C)$.
Let us prove that $\Theta$ maps $\db(\c C)$ into $\db(\c C')$. If
$X\in \db(\c C)$, then:
\begin{enumerate}
\item $\bigoplus\limits_{g\in G}Hom_{\db(\c C)}(\Sigma^m\,\t g i, X)\simeq
  Hom_{\db(A)}(\Sigma^mT_i, F_{\lambda}X)$ is finite dimensional for any
  $i\in\{1,\ldots,n\}$ and $m\in\mathbb{Z}$ (recall that
  $F_{\lambda}\colon\db(\c C)\to \db(A)$ has the covering property).
\item there exists $m_0\in \mathbb{N}$ such that $Hom_{\db(\c C)}(\Sigma^m\,\t
  g i ,X)=0$ for any $g\in G$, $i\in\{1,\ldots,n\}$ and
  $m\in\mathbb{Z}$ such that $|m|\geqslant m_0$.
\end{enumerate}
These two properties imply that $\sum\limits_{g,i,m}dim_k\ \c D(\c
C)(\Sigma^m\,\t g i , X)<\infty$. Using the fact that $(\Psi,\Theta)$ is
adjoint and using ($iv$), we deduce that $\sum\limits_{g,i,m}dim_k\ \c
D(\c C')(\Sigma^m\,_{\t g i}\c C'_?, \Theta(X))<\infty$. This proves that
$\Theta(X)\in\db(\c C')$. Therefore, $\Theta$ induces a triangle
functor $\Theta\colon\db(\c C)\to \db(\c C')$ such that the pair
$(\Psi,\Theta)$ is adjoint:
\begin{equation}
  \xymatrix{
\db(\c C')\ar@/_/@{->}[d]_{\Psi} \ar@/^/@{<-}[d]^{\Theta}\\
\db(\c C)
}\notag
\end{equation}

$\bullet$ \textbf{$\Psi\colon\db(\c C')\to\db(\c C)$ is fully faithful}.
For short, if $X,Y\in\db(\c C')$, we shall write $\Psi_{X,Y}$  for the
mapping $Hom_{\db(\c C')}(X,Y)\to Hom_{\db(\c C)}(\Psi(X),\Psi(Y))$ induced by $\Psi$.
Let $g,h\in G$ and $i,j\in\{1,\ldots,n\}$. Then $Hom_{\db(\c C')}(\
_{\t g i}\c C'_?,\ _{\t h j}\c C'_?)=Hom_{\db(\c C)}(\ \t g i,\ \t h
j)$. Moreover, we have ($iv$) $\Psi(\ _{\t g i}\c C'_?)\simeq\ \t g i$
and $\Psi(\ _{\t h j}\c C'_?)\simeq\ \t h j$, and with these
identifications, $\Psi_{_{\t g i}\c C'_?,\,_{\t h j}\c C'_?}$ is the identity mapping. On the other hand, if
$m\in\mathbb{Z}\backslash\{0\}$, then $\Psi_{\Sigma^m(\,_{\t g i}\c
  C'_?),\,_{\t h j}\c C'_?}$ is an isomorphism because
the involved morphisms spaces are trivial. Hence, $\Psi_{\Sigma^mX, Y}$ is
an isomorphism for any $m\in\mathbb{Z}$ and any projective $\c
C'$-modules $X,Y$. This shows that $\Psi$ is fully
faithful.
\null\hfill$\square$\\

\begin{lem}
\label{lem4}
  The set $\{\ \tau_{\c C}^l\Sigma^m\,\t g i\ |\ m,l\in\mathbb{Z},\
  g\in G,\ i\in\{1,\ldots,n\}\}$ generates $\db(\c C)$ as a
  triangulated category. Therefore, the functor $\Psi$ of
  Lemma~\ref{lem4} is dense.
\end{lem}
\noindent{\textbf{Proof:}} Let $\c S$ and $\sim$ be as in the proof of
Lemma~\ref{lem1}. For any $R=R_1\bigoplus\ldots\bigoplus R_n\in\c S$, fix
$\widetilde{R}_1,\ldots,\widetilde{R}_n\in\db(\c C)$ indecomposables
such that $F_{\lambda}\widetilde{R}_i\simeq R_i$ for every $i$ (see
Lemma~\ref{lem1}). Then denote by $<R>$ for the full triangulated
subcategory of $\db(\c C')$  generated by
$\{\ \tau_{\c C}^l\Sigma^m\,^g\widetilde{R}_i\ |\ m,l\in\mathbb{Z},\
  g\in G,\ i\in\{1,\ldots,n\}\}$. Remark that $<R>$ does not depend on
  the choice of $\widetilde{R}_1,\ldots,\widetilde{R}_n$ because if
  $\widetilde{R}_i'\in\db(\c C)$ verifies
  $F_{\lambda}\widetilde{R}_i'\simeq F_{\lambda}\widetilde{R}_i\simeq
  R_i$, then there exists some $g\in G$ such that
  $^g\widetilde{R}_i\simeq \widetilde{R}_i'$ (see for example the
  proof of \cite[Lem. 5.3]{lemeur6}). Since $<A>$ contains 
  all the indecomposable projective $\c C$-modules up to isomorphism,
  we infer that $<A>=\db(\c C)$. On the other hand, the second
  assertion of \cite[Prop. 6.5]{lemeur6} proves that if $R\sim R'$,
  then $<R>=<R'>$. Since $\c S$ is an equivalence class for $\sim$
  (see the proof of Lemma~\ref{lem1}), we infer that $<R>=<A>=\db(\c
  C)$ for any $R\in \c S$. This proves the first assertion of the
  lemma.

In order to prove the second assertion of the lemma, it suffices to
prove that the image of $\Psi\colon\db(\c C')\to\db(\c C)$ contains
$<T>$. By construction, this image contains $\{\,\t g i\ |\ g\in G,\
i\in\{1,\ldots,n\}\}$. On the other hand, $\Psi\colon\db(\c C')\to
\db(\c C)$ is fully faithful and admits a  right adjoint, so it preserves
Auslander-Reiten sequences. In particular, we have $\tau_{\c C}\circ
\Psi\simeq \Psi\circ \tau_{\c C'}$. This proves that the image of
$\Psi$ contains $<T>=\db(\c C)$.\hfill$\square$\\

Using Lemma~\ref{lem3} and Lemma~\ref{lem4}, the proof of
Proposition~\ref{prop1} is immediate.
Now, we are able to prove the announced
implication ($\star\star$):
\begin{prop}
\label{prop2}
  For any group $G$, there exists a connected Galois covering with
  group $G$ of $A'$ if there exists a connected Galois
  covering with group $G$ of $A'$. Consequently:
\begin{equation}
  \text{$A'=End_{\db(A)}(T)$ is simply connected}\ \Rightarrow\
  \text{$A$ is simply connected}\notag
\end{equation}
\end{prop}
\noindent{\textbf{Proof:}} Let us assume that $A'$ is simply
connected. If $F\colon \c C\to A$ is a connected
Galois covering with group $G$, then Lemma~\ref{lem2} and Proposition~\ref{prop1} show
that there exists $F'\colon\c C'\to End_{\db(A)}(T)$ a connected
Galois covering with group $G$. Since $A'$ is simply connected, we
infer (see \cite[Cor. 4.5]{lemeur2}) that $G$ is necessary the trivial group. Hence
(\textit{loc. cit.}) $A$ is simply connected.\hfill$\square$\\

\section{Hochschild cohomology and simple connectedness of squid algebras}

We refer the reader to \cite{ringel} for more details on squid
algebras. A squid algebra  over an algebraically closed field $k$
is defined by the following data: an integer $t\geqslant 2$, a
sequence $p=(p_1,\ldots,p_t)$ of non negative integers and a sequence
$\tau=(\tau_3,\ldots,\tau_t)$ of pairwise distinct non zero elements of
$k$. With this data, the squid algebra $S(t,p,\tau)$ is the
$k$-algebra $kQ/I$ where $Q$ is the following quiver:
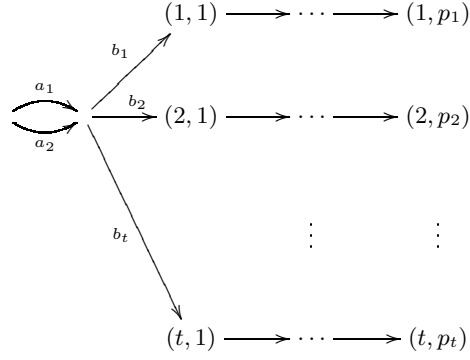
\begin{equation}
  \xymatrix{
&&(1,1)\ar@{->}[r]&\ldots \ar@{->}[r]& (1,p_1)\\
\ar@/^/@{->}[r]^{a_1}\ar@/_/@{->}[r]_{a_2} & \ar@{->}[ru]^{b_1}
\ar@{->}[r]^{b_2}
\ar@{->}[rdd]_{b_t}&(2,1)\ar@{->}[r]&\ldots \ar@{->}[r]& (2,p_2)\\
&&&\vdots&\vdots\\
&&(t,1)\ar@{->}[r]&\ldots \ar@{->}[r]& (t,p_t)
}\notag
\end{equation}

and $I$ is the ideal generated by the following relations:
\begin{equation}
  b_1a_1=b_2a_2=0,\ \ b_ia_2=\tau_i\, b_ia_1\ \ \text{for $i=3,\ldots,t$}\notag
\end{equation}

Using Happel's long exact sequence (\cite{happel}), one can compute
$HH^1(S(t,p,\tau))$:
\begin{equation}
  \operatorname{dim_k\ HH^1}(S(t,p,\tau))=\begin{cases}
1 & \text{if $t=2$}\\
0 & \text{if $t\geqslant 3$}
\end{cases}\notag
\end{equation}

On the other hand, one checks easily that if $t=2$ then the
fundamental group $\pi_1(Q,I)$ of the above
presentation of $S(t,p,\tau)$ is isomorphic to $\mathbb{Z}$ (see \cite{martinezvilla_delapena}), whereas
$S(t,p,\tau)$ is simply connected if $t\geqslant 3$. These
considerations imply the following proposition.
\begin{prop}
\label{prop3}
  Let $A$ be a squid algebra. Then $A$ is simply connected if and only
  if $HH^1(A)=0$.
\end{prop}

\section{Proof of Theorem~\ref{thm1}}

Now we can prove Theorem~\ref{thm1}. Let $A$ be quasitilted
\textit{i.e.} $A=End_{\c H}(X)^{op}$ where $\c H$ is hereditary abelian
 and where
$X\in\c H$ is basic tilting. If $\c H$ is derived equivalent to $mod(kQ)$
for some quiver $Q$, then the conclusion of the theorem follows from
\cite[Cor. 2]{lemeur6}. Otherwise, there exists $\c H'$ a
hereditary abelian category, there exists a triangle equivalence $\db(\c
H)\xrightarrow{\sim}\db(\c H')$ and there exists  $Y\in \c
H'$ basic tilting such that $End_{\c H'}(Y)^{op}$ is a squid algebra
(see \cite[Prop. 2.1, Thm. 2.6]{happel_reiten}). 
Set $B:=End_{\c H'}(Y)^{op}$. Then:
\begin{enumerate}
\item there exist triangle equivalences $\db(\c H)\xrightarrow{\sim}
  \db(A)$ and $\db(\c H')\xrightarrow{\sim}\db(B)$ mapping $X$ and $Y$
  to $A$ and $B$ respectively (thanks to \cite[Thm. 3.3, Thm 4.3]{happel_reiten_smalo}).
\item if $T\in\db(A)$ denotes the image of $Y\in\c H'$ under the equivalence
  $\db(\c H')\xleftarrow{\sim} \db(\c H)\xrightarrow{\sim}\db(A)$, then:
($i$) $Ext_A^i(T,T)=0$ for every $i\neq 0$, ($ii$) $T$ is a
cluster tilting object of $C_A$, ($iii$) $End_{\db(A)}(T)\simeq
End_{\c H'}(Y)$.
\item if $T'\in\db(B)$ denotes the image of $X\in\c H$ under the equivalence
  $\db(\c H)\xrightarrow{\sim} \db(\c H')\xrightarrow{\sim}\db(B)$, then:
($iv$) $Ext_B^i(T',T')=0$ for every $i\neq 0$, ($v$) $T'$ is a
cluster tilting object of $C_B$, ($vi$) $End_{\db(B)}(T)\simeq
End_{\c H}(X)$.
\end{enumerate}
Now, Proposition~\ref{prop2} applied  $A$ and $T$ and to $B$ and $T'$ proves
that $A$ is simply connected if and only if $B$ is simply
connected (recall that $A$ is simply connected if and only if $A^{op}$
is simply connected, see for example the proof of
\cite[Thm. 3]{lemeur6}). Since $\db(A)$ and $\db(B)$ are triangle
equivalent, we have $HH^1(A)\simeq HH^1(B)$ as $k$-vector spaces (see
\cite{keller}). Finally, Proposition~\ref{prop3} applied to $B$
proves that $A$ is simply connected if and only if $HH^1(A)=0$.\hfill$\square$\\

\bibliographystyle{plain}
\bibliography{biblio}
\end{document}